\documentstyle[amssymb,amsfonts]{amsart}

\def\R{{\hbox{\bf R}}}

\def\E{{\hbox{\bf E}}}
\def\H{{\hbox{\bf H}}}

\def\A{{\mathcal{A}}}
\def\B{{\mathcal{B}}}
\def\Energy{{\mathcal{E}}}
\font \roman = cmr10 at 10 true pt

\def\complexity{{\hbox{\roman complex}}}
\def\complexitysmall{{\hbox{\roman \scriptsize complex}}}

\def\be#1{ \begin{equation}\label{#1} }

\def\bas{\begin{align*}}
\def\eas{\end{align*}}
\def\bi{\begin{itemize}}
\def\ei{\end{itemize}}

\def\Z{{\hbox{\bf Z}}}

\def\eps{\varepsilon}
\newenvironment{proof}{\noindent {\bf Proof} }{\endprf\par}
\def \endprf{\hfill  {\vrule height6pt width6pt depth0pt}\medskip}
\def\emph#1{{\it #1}}
\def\textbf#1{{\bf #1}}

\theoremstyle{plain}
  \newtheorem{theorem}[subsection]{Theorem}

  \newtheorem{lemma}[subsection]{Lemma}
  \newtheorem{corollary}[subsection]{Corollary}

\theoremstyle{remark}
  \newtheorem{remark}[subsection]{Remark}
  
  \newtheorem{example}[subsection]{Example}

\theoremstyle{definition}
  \newtheorem{definition}[subsection]{Definition}

\include{psfig}

\begin{document}

\title[Hypergraph removal lemma]{A variant of the hypergraph removal lemma}

\author{Terence Tao}
\address{Department of Mathematics, UCLA, Los Angeles CA 90095-1555}
\email{tao@@math.ucla.edu}

\begin{abstract}  Recent work of Gowers \cite{gowers} and Nagle, R\"odl, Schacht, and Skokan \cite{nrs}, \cite{rodl}, \cite{rodl2} has established
a hypergraph removal lemma, which in turn implies some results of Szemer\'edi \cite{szemeredi} and Furstenberg-Katznelson \cite{fk} concerning one-dimensional and multi-dimensional arithmetic progressions respectively.  In this paper we shall give a self-contained proof of this hypergraph removal lemma.  In fact we prove a slight strengthening of the result, which we will use in a subsequent paper \cite{tao-multiprime} to establish (among other things) infinitely many constellations of a prescribed shape in the Gaussian primes.   
\end{abstract}

\maketitle

\section{Introduction}

In this paper we prove a slight variant of the hypergraph removal lemma established recently and independently by Gowers \cite{gowers}
and Nagle, R\"odl, Schacht and Skokan \cite{nrs}, \cite{rodl}, \cite{rodl2}.  To motivate this lemma, let us first recall the more well-known triangle removal lemma from graph theory of
Ruzsa and Szemer\'edi \cite{rsz}.  It will be convenient to work in the setting of tripartite graphs, though we will comment about the generalization
to general graphs shortly.  We adopt the following $o()$ and $O()$ notation:  If $x,y_1,\ldots,y_n$ are parameters,
we use $o_{x \to 0; y_1,\ldots,y_n}(X)$ to denote any quantity bounded in magnitude by $X c(x,y_1,\ldots,y_n)$, where $c()$ is a function which goes to zero as $x \to 0$ for each fixed choice of $y_1,\ldots,y_n$.  Similarly, we use $O_{y_1,\ldots,y_n}(X)$ to denote any quantity bounded by
$X C(y_1,\ldots,y_n)$, for some function $C()$ of $y_1,\ldots,y_n$.  If $A$ is a finite set, we use $|A|$ to denote the cardinality of $A$.

\begin{theorem}[Triangle removal lemma, tripartite graph version]\label{triangle-removal}\cite{rsz}  Let $V_1, V_2, V_3$ be finite non-empty sets of vertices, and let $G = (V_1,V_2,V_3,E_{12}, E_{23}, E_{31})$ be a tri-partite graph on these sets of vertices, thus $E_{ij} \subseteq V_i \times V_j$ for $ij = 12,23,31$.  Suppose that the number of triangles in this graph does not exceed $\delta |V_1| |V_2| |V_3|$ for some $0 < \delta < 1$.
Then there exists a graph $G' = G'(V_1,V_2,V_3,E'_{12}, E'_{23}, E'_{31})$ which contains no triangles whatsoever, and 
such that $|E_{ij} \backslash E'_{ij}|= o_{\delta \to 0}(|V_i \times V_j|)$ for $ij=12,23,31$.
\end{theorem}

One can view $G'$ as a ``triangle-free approximation'' to $G$.
Note that we do not assume that $G'$ is a subgraph of $G$, but one can easily obtain this conclusion by replacing $E'_{ij}$ with $E'_{ij} \cap E_{ij}$ if desired (i.e. one replaces $G'$ by $G' \cap G$).  As we shall see, however, it will be convenient to allow the 
possibility that $G'$ is not a subgraph of $G$.

\begin{remark}
The above theorem is phrased for tri-partite graphs, but it quickly implies an analogous version for non-partite graphs $G = (V,E)$, by
taking three copies $V_1 = V_2 = V_3 = V$ of the vertex set $V$, and constructing the
bipartite graph $\tilde G = (V_1,V_2,V_3,E_{12},E_{23},E_{31})$, where $E_{ij}$ consists of those pairs $( x, y )$ which
are the endpoints of an edge in $E$.  We omit the details.
\end{remark}

It was observed in \cite{rsz} that Theorem \ref{triangle-removal} implies Roth's famous theorem \cite{roth} that subsets of integers of positive density
contain infinitely many progressions of length three.  In \cite{soly-roth} it was also observed that Theorem \ref{triangle-removal} also implies
that subsets of $\Z^2$ with positive density contain infinitely many right-angled triangles (a result
first obtained in \cite{AS}).  It was observed earlier (for instance in \cite{rodl-icm} or \cite{frankl02})
that an extension of the triangle removal lemma to hypergraphs would similarly imply
Szemer\'edi's famous theorem \cite{szemeredi} on progressions of arbitrary length; by modifying the observation in \cite{soly-roth}, it
would also imply a multidimensional extension of that theorem due to Furstenberg and Katznelson \cite{fk}.  We shall return to this issue
in the sequel \cite{tao-multiprime} to this paper, and discuss the above hypergraph removal lemma in detail later in this introduction.

Theorem \ref{triangle-removal} was proven using the \emph{Szemer\'edi regularity lemma} (see e.g. \cite{szemeredi-reg}, \cite{komlos} for a survey of this lemma and its applications), which roughly speaking allows one to approximate an arbitrary large and complex graph to arbitrary accuracy by a much simpler object; see also \cite{van}, \cite{shkredov} for further refinements of Theorem \ref{triangle-removal}.  This proof in fact yields a little bit
more information on the triangle-free approximation $G'$ to $G$, namely that $G'$ can be chosen to be ``bounded complexity''.  More precisely:

\begin{theorem}[Strong triangle removal lemma, tripartite graph version]\label{triangle-removal-2}\cite{rsz}  Let $V_1, V_2, V_3$ be finite non-empty sets of vertices, and let $G = (V_1,V_2,V_3,E_{12}, E_{23}, E_{31})$ be a tri-partite graph on these sets of vertices.  Suppose that $G$ contains
at most $\delta |V_1| |V_2| |V_3|$ triangles.
Then there exists a graph $G' = G'(V_1,V_2,V_3,E'_{12}, E'_{23}, E'_{31})$ which contains no triangles whatsoever, and 
such that $|E_{ij} \backslash E'_{ij}|= o_{\delta \to 0}(|V_i \times V_j|)$ for $ij=12,23,31$.  Furthermore, there exists a quantity
$M = O_\delta(1)$, and partitions $V_i = V_{i,1} \cup \ldots V_{i,M}$ for each $i=1,2,3$ into sets $V_{i,a}$ (some of which may be empty)
such that for each $ij=12,23,31$, $E'_{ij}$ is the union of sets of the form $V_{i,a} \times V_{j,b}$.
\end{theorem}

Note that the graph $G'$ constructed in Theorem \ref{triangle-removal-2} will typically not be a subgraph of $G$.  One could make the sets
$V_{i,1},\ldots,V_{i,M}$ to be the same size (with at most one exception for each $i$) without much difficulty but we will not endeavour to do so here.
There is also a version of this lemma for non-tripartite graphs which is well known (and essentially equivalent to the tripartite version)
but we will not reproduce it here.

It turns out that Theorem \ref{triangle-removal} and Theorem \ref{triangle-removal-2} can be rephrased in a more ``probabilistic''
manner.  One reason for doing this is because in our arguments we will need two basic concepts
from probability theory, which are  \emph{conditional expectation} and \emph{complexity} respectively.  It
seems that with the aid of these concepts, the proofs become somewhat cleaner to give\footnote{For a more traditional combinatorial approach to these problems, see \cite{rs}.}.  
To explain these concepts we need some notation.  For reasons which will become clearer later, we shall use a rather general notation which
incorporates the above Theorems as a special case.

\begin{definition}[Hypergraphs] If $J$ is a finite set and $d \geq 0$, we define ${J \choose d} := \{ e \subseteq J: |e| = d \}$ to be the set of all
subsets of $J$ of cardinality $d$.  A \emph{$d$-uniform hypergraph} on $J$ is then defined to be any subset $H_d \subseteq {J \choose d}$ of
${J \choose d}$.  For instance, an undirected graph $G = (V,E)$ without loops can be viewed as a $2$-uniform hypergraph on $V$.
\end{definition}

\begin{example}\label{triangle-ex} If $J := \{1,2,3\}$, then the triangle $H_2 := {J \choose 2} = 
\{\{1,2\}, \{2,3\}, \{3,1\}\}$ is a 2-uniform hypergraph on $J$.
\end{example}

\begin{definition}[Hypergraph systems]  A \emph{hypergraph system} is a quadruplet $V = (J, (V_j)_{j \in J}, d, H_d)$, where $J$ is a finite set,
$(V_j)_{j \in J}$ is a collection of finite non-empty sets indexed by $J$, $d \geq 1$ is positive integer, and $H_d \subseteq {J \choose d}$ is a
$d$-uniform hypergraph.  For any $e \subseteq J$, we set $V_e := \prod_{j \in e} V_j$, and let $\pi_e: V_J \to V_e$ be the canonical projection map.
\end{definition}

\begin{remark}
Very roughly speaking, a hypergraph system corresponds to the notion of a \emph{measure-preserving system}\footnote{A measure preserving system is a probability space $(X, {\mathcal B}, \mu)$ together with a shift $T: X \to X$ that preserves the measure $\mu$.  The ergodic approach to Szemer\'edi's theorem, as introduced by Furstenberg\cite{furst}, recasts the problem of finding arithmetic progressions as that of understanding averages such as $\liminf_{N \to \infty} \frac{1}{N} \sum_{n=1}^N \mu(A \cap T^n A \cap \ldots \cap T^{(k-1)n} A)$.  This can in turn be viewed as the problem of understanding shift operators such as $(T, T^2,\ldots,T^{k-1})$ on a product space $X \times \ldots \times X$.  This has some intriguing parallels with the combinatorial approach, in which the problem of obtaining arithmetic progressions in a set $V$ is reduced to that of analyzing Cayley-type graphs or hypergraphs, which can be viewed as subsets of $V \times \ldots \times V$.  We do not know of any formal connection between these two approaches, nevertheless there do appear to be some interesting similarities.}
 in ergodic theory, though with the notable difference that no analogue of the shift operator exists in a hypergraph system.  Indeed the $V_j$ are simply finite sets, and need not have any additive structure whatsoever.  
\end{remark}

\begin{definition}[Conditional expectation]  Let $V = (J, (V_j)_{j \in J}, d, H_d)$ be a hypergraph system.  
If $f: V_J \to \R$ is a function, we define the expectation $\E(f) = \E(f(x) | x \in V_J)$ by the formula
$$ \E(f) = \E(f(x) | x \in V_J) := \frac{1}{|V_J|} \sum_{x \in V_J} f(x).$$
Similarly, if $\B$ is a $\sigma$-algebra\footnote{Of course, since $V_J$ is finite, we do not need to distinguish finite unions and countable unions, and could simply call $\B$ an ``algebra'', or even a ``partition''; the latter notation is in fact used in most treatments of the regularity lemma.  However we prefer the notation of $\sigma$-algebra as being highly suggestive, evoking ideas and insights from probability theory, measure theory, and information theory.}
 on $V_J$, i.e. a collection of sets in $V_J$ which contains $\emptyset$ and $V_J$, and
is closed under unions, intersections, and complementation, we define the \emph{conditional expectation} $\E(f|\B): V_J \to \R$ by the formula
$$ \E(f|\B)(x) := \frac{1}{|\B(x)|} \sum_{y \in \B(x)} f(y),$$
where $\B(x)$ is the smallest element of $\B$ which contains $x$.  For each $e \subseteq J$, let $\A_e$ be the $\sigma$-algebra on $V_J$ 
defined by $\A_e := \{ \pi_e^{-1}(E): E \subseteq V_e \}$.  In other words, $\A_e$ consists of those subsets of $V_J$, 
membership of which is determined solely by the co-ordinates of $V_J$ indexed by $e$.
\end{definition}

One can interpret the usage of these averages as imposing the uniform probability distribution on each $V_e$, which basically amounts
to introducing a set $(x_j)_{j \in J}$ of independent random variables, with each $x_j$ ranging uniformly in $V_j$.  

If $\B_1$ and $\B_2$ are two $\sigma$-algebras on $V_J$, we use $\B_1 \vee \B_2$ to denote the smallest $\sigma$-algebra that contains both $\B_1$ and $\B_2$; this corresponds to the familiar concept of the \emph{common refinement} of two partitions.  We can more generally define $\bigvee_{i \in I} \B_i$ for any collection $(\B_i)_{i \in I}$ of $\sigma$-algebras.

\begin{example} For any finite non-empty sets $V_1,V_2,V_3$, the quadruplet $V = (J, (V_j)_{j \in J}, 2, H_2)$ is a hypergraph system,
where $J := \{1,2,3\}$ and $H_2 := {J \choose 2}$ are as in Example \ref{triangle-ex}.  
The $\sigma$-algebra $\A_{\{1,2\}}$ is the algebra of all subsets of $V_1 \times V_2 \times V_3$
which do not depend on the third variable, and thus take the form $E \times V_3$ for some $E \subseteq V_1 \times V_2$.  Similarly for $\A_{\{2,3\}}$
and $\A_{\{3,1\}}$.
\end{example}

\begin{definition}[Complexity]  Let $V = (J, (V_j)_{j \in J}, d, H_d)$ be a hypergraph system.
If $\B$ is a $\sigma$-algebra in $V_J$, we define the \emph{complexity} $\complexity(\B)$ of $\B$ to be the least number of sets in $V_J$ 
needed to generate $\B$ as a $\sigma$-algebra; this can be viewed as a simplified version of the Shannon entropy $\H(\B)$, which we will not use here.  We observe the obvious inequalities
\begin{equation}\label{complex-jump}
\complexity(\B_1 \vee \B_2) \leq \complexity(\B_1) + \complexity(\B_2) \hbox{ for arbitrary } \B_1, \B_2
\end{equation}
and
\begin{equation}\label{b-card}
 |\B| \leq 2^{2^{\complexitysmall(\B)}}.
\end{equation}
\end{definition}

\begin{remark} If one views $\B$ as a partition, the complexity is essentially the logarithm of the number of cells in the partition.
From an information-theoretic perspective, the complexity measures how many bits of information are needed to know which atom of $\B$ a given point in $V_J$ lies in.
\end{remark}

If $E$ is a subset of $V_J$, we let $1_E: V_J \to \R$ be the indicator function, thus $1_E(x) := 1$ when $x\in E$ and $1_E(x) := 0$ otherwise.
In particular, $\E(1_E) = |E|/|V_J|$ can be viewed as the ``density'' or ``probability'' of $E$ in $V_J$.

With all this notation, Theorem \ref{triangle-removal-2} becomes

\begin{theorem}[Strong triangle removal lemma, $\sigma$-algebra version]\label{triangle-main}  Let $V = (J, (V_j)_{j \in J}, d, H_d)$ be
a hypergraph system with $J = \{1,2,3\}$, $d = 2$, and $H_d = {J \choose d} = \{ \{1,2\}, \{2,3\}, \{3,1\} \}$.
For each $e \in H_d$, let $E_e$
be a set in $\A_e$ such that
$$ \E( \prod_{e \in H_d} 1_{E_e} ) \leq \delta$$  
for some $0 < \delta < 1$.  Then there exist sets $E'_e \in \A_e$ for $e \in H_d$
such that
$$ \bigcap_{e \in H_d} E'_e = \emptyset$$
and 
$$\E( 1_{E_e \backslash E'_e} ) = o_{\delta \to 0}(1) \hbox{ for all } e \in H_d.$$
Furthermore, for each $i \in J$ there exists sub-algebras $\B_i \subseteq \A_{\{i\}}$ 
such that
$$ \complexity(\B_i) = O_\delta(1) \hbox{ for } i \in J$$
and 
$$ E'_e \in \bigvee_{i \in e} \B_i \hbox{ for } e \in H_d.$$
\end{theorem}

It is easy to see that Theorem \ref{triangle-removal-2} and Theorem \ref{triangle-main} are equivalent.  The notation here may appear
quite cumbersome, but the advantages of these notations will hopefully become more apparent when we prove a generalization of this result shortly.

The case of $d=2$, and $J$ and $H_d$ arbitrary, was treated in \cite{efr}.  It was then conjectured in that paper that a result of the above
type should also hold for higher $d$.  The generalization of Theorem \ref{triangle-removal} to the higher $d$ case 
was accomplished only recently and independently by Gowers \cite{gowers-hyper} and Nagle, R\"odl, Schacht, Skokan \cite{nrs}, \cite{rodl}, \cite{rodl2},
using the language of hypergraphs.  It turns out that Theorem \ref{triangle-removal-2} or Theorem \ref{triangle-main} can similarly be generalized, and with the notation already developed, the extension is very easy to state:

\begin{theorem}[Hypergraph removal lemma]\label{main-2}\cite{gowers-hyper}, \cite{nrs}, \cite{rodl}, \cite{rodl2}  
Let $V = (J, (V_j)_{j \in J}, d, H_d)$ be a hypergraph system.  For each $e \in H$, let $E_e$ be a set in $\A_e$ such that
\begin{equation}\label{E-dens}
 \E( \prod_{e \in H_d} 1_{E_e} ) \leq \delta
 \end{equation}
for some $0 < \delta < 1$.  Then for each $e \in H_d$ there exists a set $E'_e \in \A_e$ such that
\begin{equation}\label{E-cap}
 \bigcap_{e \in H_d} E'_e = \emptyset
 \end{equation}
and
\begin{equation}\label{E-error}
\E( 1_{E_e \backslash E'_e} ) = o_{\delta \to 0; J}(1) \hbox{ for all } e \in H_d.
\end{equation}
Furthermore, there exist sub-algebras $\B_{e'} \subseteq \A_{e'}$ whenever $e' \subset J$ and $|e'| < d$ obeying the complexity estimate
\begin{equation}\label{E-complex}
\complexity(\B_{e'}) = O_{J, \delta}(1) \hbox{ whenever } e' \subseteq J \hbox{ and } |e'| < d
\end{equation}
(so in particular $|\B_{e'}| = O_{J, \delta}(1)$, thanks to \eqref{b-card})
and
\begin{equation}\label{E-meas}
E'_e \in \bigvee_{e' \subsetneq e} \B_{e'} \hbox{ for all } e \in H_d.
\end{equation}
\end{theorem}

Clearly Theorem \ref{triangle-main} is a special case of Theorem \ref{main-2}.
We have attributed this theorem to Gowers \cite{gowers-hyper} and Nagle-R\"odl-Schacht-Skokan \cite{nrs}, \cite{rodl}, \cite{rodl2} because it 
follows from their methods, although a theorem of this type is not stated explicitly in those papers.  One can formulate variants of this 
removal lemma in the case when $H_d$ is not $d$-uniform but we will not do so here.  A related result has recently been obtained
in \cite{rs}, using techniques similar in spirit to those here (though with substantially different notation).

The main purpose of this paper is 
to explicitly prove Theorem \ref{main-2} in a completely self-contained manner. 
In a subsequent paper \cite{tao-multiprime}, we will then transfer this theorem (as in \cite{gt-primes}) to
obtain a relative version of Theorem \ref{main-2}, restricted to a suitably pseudorandom subset of $\prod_j V_j$.  This will
then be used (again following \cite{gt-primes}) to deduce 
the existence of infinitely many constellations of a prescribed shape in 
the Gaussian primes and similar sets.  

As a corollary of Theorem \ref{main-2}, we obtain the hypergraph removal lemma in a formulation closer to that of Gowers or Nagle-R\"odl-Schacht-Skokan:

\begin{corollary}[Hypergraph removal lemma, partite hypergraph version]\label{hypergraph-removal}\cite{gowers-hyper}, \cite{nrs},\cite{rodl}, \cite{rodl2}  Let $(V_j)_{j \in J}$ be a collection of finite non-empty sets. Let $0 \leq d \leq |J|$, and let $H_d \subseteq {J \choose d}$ 
be a $d$-uniform hypergraph on $J$.  
For each $e \in H_d$, let $E_e$ be a subset of $\prod_{j \in e} V_j$.  Suppose that
$$ |\{ (x_j)_{j \in J} \in \prod_{j \in J} V_j: (x_j)_{j\in e} \in E_e \hbox{ for all } e \in H_d \}|
\leq \delta \prod_{j \in J} |V_j|$$
for some $0 < \delta \leq 1$; in other words, the $J$-partite hypergraph $G = ((V_j)_{j \in J},(E_e)_{e \in H_d})$ contains at most
$\delta \prod_{j \in J} |V_j|$ copies of $H_d$.  Then for each $e \in H_d$ there exists $E'_e \subset \prod_{j \in e} V_j$ such that
$$ \{ (x_j)_{j \in J} \in \prod_{j \in J} V_j: (x_j)_{j\in e} \in E'_e \hbox{ for all } e \in H_d \} = \emptyset$$
(i.e. the $J$-partite hypergraph $G' = G'((V_j)_{j \in J},(E'_e)_{e \in H_d})$ contains no copies of $H_d$ whatsoever), and
such that
$|E_e \backslash E'_e| = o_{\delta \to 0; |J|}(\prod_{j \in e} |V_j| )$ for all $e \in H_d$.
\end{corollary}

The deduction of Corollary \ref{hypergraph-removal} from Theorem \ref{main-2} is analogous to the deduction of Theorem \ref{triangle-removal}
from Theorem \ref{triangle-main} and is omitted.  It seems quite likely that we can obtain similar analogues for non-partite hypergraphs,
just as was the case with the non-partite version of Theorem \ref{triangle-removal}; see \cite{gowers-hyper}, \cite{nrs}, \cite{rodl}, \cite{rodl2} for
some examples of this, though for applications to Szemer\'edi-type theorems it is the partite version which is of importance.  
It should be unsurprising that
Theorem \ref{triangle-removal} is then the special case of Corollary \ref{hypergraph-removal} applied to the (hyper)graph 
in Example \ref{triangle-ex}.  
The case $|J|=4$ and $H_3 = {J \choose 3}$  was
treated in \cite{frankl02}. Just as Theorem \ref{triangle-removal} implies Roth's theorem, Corollary \ref{hypergraph-removal} implies Szemer\'edi's theorem \cite{szemeredi} on arithmetic progressions, as well as the multidimensional generalization of that theorem due to Furstenberg and Katznelson \cite{fk}; see \cite{soly-2}, \cite{frankl02}, \cite{gowers-hyper}, \cite{rodl2} for further discussion\footnote{It was also recently observed that this hypergraph removal result also implies another theorem of Furstenberg and Katznelson \cite{fk2} on affine subspaces of dense subsets of high-dimensional finite field vector spaces; see \cite{rstt}.}.  Thus this paper provides
a moderately short and self-contained proof of these theorems, although we emphasize that this goal was already achieved in the prior work
of \cite{gowers-hyper}, \cite{nrs}, \cite{rodl}, \cite{rodl2}.

The remainder of this paper is devoted to proving Theorem \ref{main-2}.  As one might expect from
the previous proofs of these types of results, our proof shall proceed by proving a ``hypergraph regularity lemma'' and a ``hypergraph counting
lemma''.  The arguments are broadly along similar lines to those of Gowers or Nagle, R\"odl, Schacht, and Skokan, although it seems that using the notation of
$\sigma$-algebras and probability theory allows for slightly cleaner arguments.

The author thanks Fan Chung Graham, Vojt\v{e}ch R\"odl, Mathias Schacht, and Jozsef Solymosi
for helpful comments and references.  He is particularly indebted to
 Mathias Schacht for supplying the recent preprint \cite{rs}, and to the anonymous referees for a careful reading of the paper and
 many cogent suggestions and corrections.  The author is supported by a grant from the Packard foundation.

\section{Pseudorandomness and the regularity lemma}

Henceforth the hypergraph system $V = (J, (V_j)_{j \in J}, d, H_d)$ will be fixed.  In this section we shall state and prove
a $\sigma$-algebra version of the hypergraph regularity lemma (Lemma \ref{full-regularity}).  
This lemma establishes a dichotomy between pseudorandomness (or $\eps$-regularity, or small discrepancy) on one hand, and bounded complexity\footnote{This is very similar to the dichotomy between weak mixing and compactness in ergodic theory, which is of great utility in proving statements such as Szemer\'edi's theorem; it seems of interest to explore these connections further.} on the other; the regularity lemma then asserts, very roughly speaking,
that any given set or $\sigma$-algebra (or family of $\sigma$-algebras) can be split into a component with bounded complexity, and a 
component which is pseudorandom (has small discrepancy).  

In order to state the regularity lemma we need to formalize the notion of pseudorandomness (or more precisely, of discrepancy).  We shall
also need a notion of the \emph{energy} of a $\sigma$-algebra in order to keep track of the inductions that go into the
proof of the regularity lemma, and also in the final statement of our regularity lemma.

We shall not state the final regularity lemma we need (Lemma \ref{full-regularity}) immediately.  To begin with, we set out
our notation for discrepancy and energy.  Initially we shall be focusing primarily on a single edge $e \subseteq J$, as opposed to an entire
hypergraph $H_d$, though this hypergraph shall emerge later in this section.

\begin{definition}[$e$-discrepancy]
For any $e \subseteq J$, we define the \emph{skeleton} $\partial e$ of $e$ to be the set $\{ f \subsetneq e: |f| = |e|-1\}$.
If $e \subseteq J$, $E_e \subseteq V_J$, and $\B$ is a $\sigma$-algebra on $V_J$, 
we define the \emph{$e$-discrepancy} $\Delta_e(E_e|\B)$ of the set $E_e$ with respect to the $\sigma$-algebra $\B$ to be
the quantity\footnote{This quantity is related to the Gowers uniformity norms used for instance in \cite{gowers}, \cite{gowers-hyper}, \cite{gt-primes}, but we will not explicitly introduce those norms here.  This quantity is also related to the notion of a pseudorandom hypergraph, studied for instance in \cite{krs-hyper}.}
\begin{equation}\label{Psie}
\Delta_e(E_e|\B) := \sup_{E_f \in \A_f \forall f \in \partial e} 
| \E\left( (1_{E_e} - \E(1_{E_e}|\B)\right) \prod_{f \in \partial e} 1_{E_f} )|
\end{equation}
where the supremum is over all collections of sets $(E_f)_{f \in \partial e}$, where each $E_f$ lies in the $\sigma$-algebra $\A_f$.
Note that since $V_J$ is finite, so is $\Delta_e(E_e|\B)$.
\end{definition}

Roughly speaking, the $e$-discrepancy $\Delta_e(E_e|\B)$ measures the amount of ``structure'' in $E_e$ which is not already captured by the $\sigma$-algebra
$\B$.  By ``structure'', we mean sets which can be easily described by sets from the lower order $\sigma$-algebras $\A_f$, as opposed to a
generic set in $\A_e$ which in general is likely to have no good decomposition (or approximate decomposition) into sets from the $\A_f$.
Thus if $\Delta_e(E_e|\B)$ is small, we expect $E_e$ to behave randomly (i.e. in an unstructured way) on most atoms of $\B$. 
The $\Delta_e(E_e|\B)$
generalize the concept of $\eps$-regularity, as the following
example shows:

\begin{example}\label{gve}  Let $G = (V_1, V_2, E_{12})$ be a bipartite graph between two finite non-empty sets $V_1, V_2$; we can thus view
$E_{12}$ as a set in $\A_{\{1,2\}}$, 
where $V$ is the hypergraph system $V = (J, (V_j)_{j \in J}, d, H_d)$ with $J = \{1,2\}$, $d=2$, and $H_d = {J \choose d} =
\{ \{1,2\} \}$.    Suppose that
$E_{12}$ has density $\E(1_{E_{12}}) = \sigma$ (i.e. $\sigma = |E_{12}|/|V_1| |V_2|$), and that
$$ \Delta_{\{1,2\}}(E_{12}|\A_\emptyset) \leq \eps$$
for some $\eps > 0$.  Then by definition we have
$$ |\E( (1_{E_{12}} - \sigma) 1_{E_1} 1_{E_2} )| \leq \eps \hbox{ whenever } E_1 \in \A_{\{1\}}, E_2 \in \A_{\{2\}}.$$
In the original setting of the bipartite graph $G$, this is equivalent to asserting that
$$ \bigl| |E_{12} \cap (E_1 \times E_2)| - \sigma |E_1| |E_2| \bigr| \leq \eps |V_1| |V_2|$$
for all $E_1 \subseteq V_1$ and $E_2 \subseteq V_2$.  The reader may recognize this as a pseudorandomness condition or
$\eps$-regularity condition on the graph $G$.  If we replace $\A_\emptyset$ by a finer $\sigma$-algebra such as $\B_1 \vee \B_2$ for some
$\B_1 \subseteq \A_{\{1\}}$ and $\B_2 \subseteq \A_{\{2\}}$, where the complexity of $\B_1$ and $\B_2$ is small compared to $1/\eps$,
then a condition such as $\Delta_{\{1,2\}}(E_{12}|\B_1 \vee \B_2) \leq \eps$ states, roughly speaking, that the graph $G$ is $\eps$-regular on ``most''
of the atoms $A_1 \times A_2$ in the partition associated to $\B_1 \vee \B_2$.
\end{example}

If $\B$ is a $\sigma$-algebra on $V_J$ and $E$ is a set in $V_J$ (not necessarily in $\B$), 
we define the \emph{$E$-energy} of $\B$ to be the quantity
$$ \Energy_E(\B) := \E( |\E(1_E|\B)|^2 ).$$
Clearly, the $E$-energy $\Energy_E(\B)$ ranges between 0 and $1$; intuitively, $\Energy_E(\B)$ is a measure of how much information
about $E$ is captured by $\B$, and is thus in many ways complementary to the $e$-discrepancy $\Delta_e(E|\B)$.  From Pythagoras' theorem we can verify the identity
\begin{equation}\label{pythagoras}
\Energy_E(\B') = \Energy_E(\B) + \E( |\E(1_E|\B') - \E(1_E|\B)|^2 ) \hbox{ whenever } \B \subseteq \B',
\end{equation}
thus finer $\sigma$-algebras have larger $E$-energy.

\begin{remark} In the setting of Example \ref{gve} with $\B = \B_1 \vee \B_2$ for some $\B_1 \subseteq \A_{\{1\}}$
and $\B_2 \subseteq \A_{\{2\}}$, the energy is a familiar quantity in the theory of the regularity lemma, and is usually referred
to as the \emph{index} of the partition; see \cite{szemeredi-reg}.
\end{remark}

Let us informally say that a set $E_e \in \A_e$ is \emph{$e$-pseudorandom with respect to $\B$} if the $e$-discrepancy $\Delta_e(E_e|\B)$ is small.
A fundamental fact (which was already exploited in \cite{szemeredi}, \cite{szemeredi-reg}) 
is that if $E$ is \emph{not} $e$-pseudorandom with respect to $\B$, then we can find a refinement of $\B$ with 
higher energy and not much larger complexity:

\begin{lemma}[Large discrepancy implies energy increment]\label{increment}  
Let $e \subseteq J$, let $E_e \in \A_e$ be a set, and for each $f \in \partial e$ let
$\B_f \subseteq \A_f$ be a $\sigma$-algebra such that
$$ \Delta_e(E_e|\bigvee_{f \in \partial e} \B_f) \geq \eps$$
for some $\eps > 0$.  Then there exists a $\sigma$-algebra $\B_f \subseteq \B'_f \subseteq \A_f$ for all $f \in \partial e$ such that
\begin{equation}\label{complexity-double}
 \complexity( \B'_f ) \leq \complexity( \B_f ) + 1
 \end{equation}
and
\begin{equation}\label{energy-increment}
 \Energy_{E_e}(\bigvee_{f \in \partial e} \B'_f) \geq \Energy_{E_e}(\bigvee_{f \in \partial e} \B_f) + \eps^2.
 \end{equation}
\end{lemma}

\begin{proof}  By \eqref{Psie} (and the finiteness of $V_J$) we can find sets $E_f \in \A_f$ for all $f \in \partial e$ such that
$$ |\E\left( \bigl(1_{E_e} - \E(1_{E_e}|\bigvee_{f \in \partial e} \B_f)\bigr) \prod_{f \in \partial e} 1_{E_f} \right)| \geq \eps.$$
For each $f \in \partial e$, let $\B'_f$ be the $\sigma$-algebra
$$ \B'_f := \B_f \vee \B(E_f)$$
then we have $\B_f \subseteq \B'_f \subseteq \A_f$, and obtain \eqref{complexity-double} from \eqref{complex-jump}.
Since $\prod_{f \in \partial e} 1_{E_f}$ is measurable with respect to $\bigvee_{f \in \partial e} \B'_f$,
and $1_{E_e} - \E(1_{E_e}|\bigvee_{f \in \partial e} \B'_f)$ has zero conditional expectation with respect to $\bigvee_{f \in \partial e} \B'_f$ we see that
$$ \E\left( \bigl(1_{E_e} - \E(1_{E_e}|\bigvee_{f \in \partial e} \B'_f)\bigr) \prod_{f \in \partial e} 1_{E_f} \right) = 0$$
and hence
$$ |\E\left( \bigl(\E(1_{E_e}|\bigvee_{f \in \partial e} \B'_f) - \E(1_{E_e}|\bigvee_{f \in \partial e} \B_f)\bigr) 
\prod_{f \in \partial e} 1_{E_f} \right)| \geq \eps.$$
By the boundedness of $\prod_{f \in \partial e} 1_{E_f}$ and the Cauchy-Schwarz inequality we conclude
$$ \E\left( \bigl|\E(1_{E_e}|\bigvee_{f \in \partial e} \B'_f) - \E(1_{E_e}|\bigvee_{f \in \partial e} \B_f)\bigr|^2 \right) \geq \eps^2,$$
and \eqref{energy-increment} then follows from \eqref{pythagoras}.
\end{proof}

By iterating Lemma \ref{increment}, one expects to be able to show that any given set $E_e \in \A_e$ must be $e$-pseudorandom with respect to a $\sigma$-algebra $\B$ of bounded complexity, since otherwise we could create a tower of $\sigma$-algebras whose energy increments indefinitely.  
Such statements can be viewed as $\sigma$-algebra analogues of the Szemer\'edi regularity lemma.  There are several such lemmas available; the
final lemma which we need is a bit lengthy to state, so we begin by stating some simpler regularity lemmas which we will then iterate to obtain
the stronger lemmas which we need.  We first obtain a preliminary iteration of Lemma \ref{increment}, in which the 
single set $E_e \in A_e$ is replaced by an ensemble of sets, or more precisely an 
ensemble $(\B_e)_{e \in H}$ of $\sigma$-algebras with bounded complexity.

If $H_d$ is a $d$-uniform hypergraph, we define $\partial H_d$ to be the $(d-1)$-uniform hypergraph
$\partial H_d := \bigcup_{e \in H_d} \partial e$.  

\begin{lemma}[Dichotomy between randomness and structure]\label{dichotomy}  Let $V = (J, (V_j)_{j \in J}, d, H_d)$ be a hypergraph system.
For each $e \in H_d$, let $\B_e \subseteq \A_e$ be a $\sigma$-algebra with the complexity bounds
$$ \complexity(\B_e) \leq m \hbox{ for all } e \in H_d$$
for some $m > 0$, and for each $f \in \partial H_d$, let $\B_f \subseteq \A_f$ be a $\sigma$-algebra with the complexity bounds
$$ \complexity(\B_f) \leq M \hbox{ for all } f \in \partial H_d$$
for some $M > 0$.  Let $\eps, \delta > 0$.  Then one of the following statements must hold.

\begin{itemize}
\item (Randomness) There exists $\sigma$-algebras $\B_f \subseteq \B'_f \subseteq \A_f$ for all $f \in \partial H_d$ such that
\begin{equation}\label{ebe-1}
\Energy_{E_e}(\bigvee_{f \in \partial e} \B'_f) < \Energy_{E_e}(\bigvee_{f \in \partial e} \B_f) + \eps^2 \hbox{ for all } e \in H_d \hbox{ and } E_e \in \B_e
\end{equation}
and
\begin{equation}\label{ebe-2}
\Delta_e(E_e|\bigvee_{f \in \partial e} \B'_f) \leq \delta \hbox{ for all } e \in H_d \hbox{ and } E_e \in \B_e.
\end{equation}

\item (Structure) There exist $\sigma$-algebras $\B_f \subseteq \B'_f \subseteq \A_f$ for all $f \in \partial H_d$ such that
\begin{equation}\label{ebe-3}
\Energy_{E_e}(\bigvee_{f \in \partial e} \B'_f) \geq \Energy_{E_e}(\bigvee_{f \in \partial e} \B_f) + \eps^2 \hbox{ for some } e \in H_d \hbox{ and } E_e \in \B_e
\end{equation}
and
\begin{equation}\label{ebe-4}
\complexity(\B'_f) \leq M + O_{|J|, m, \eps, \delta}(1) \hbox{ for all } f \in \partial H_d.
\end{equation}
\end{itemize}
\end{lemma}

\begin{proof}  We run the following algorithm:

\begin{itemize}
\item Step 0.  Initialize $\B'_f := \B_f$ for all $f \in \partial H_d$.  Note that \eqref{ebe-1} and \eqref{ebe-4} currently hold.
\item Step 1.  If \eqref{ebe-2} holds, then we halt the algorithm (we are in the ``randomness'' half of the dichotomy).  
Otherwise, there exists an $e \in H$ and $E_e \in \B_e$ such that
$$ \Delta_e(E_e|\bigvee_{f \in \partial e} \B'_f) > \delta.$$
We can then invoke Lemma \ref{increment} to locate 
refinements $\B'_f \subseteq \B''_f \subseteq \A_f$ for all $f \in \partial H_d$ (note that $\B''_f$ will just equal $\B'_f$ if
$f \not \subset e$) such that
$$ \complexity(\B''_f) \leq \complexity(\B'_f) + 1 \hbox{ for all } f \in \partial H_d$$
and
$$ \Energy_{E_e}(\bigvee_{f \in \partial e} \B''_f) \geq \Energy_{E_e}(\bigvee_{f \in \partial e} \B'_f) + \delta^2.$$
\item Step 2.  We replace $\B'_f$ with $\B''_f$ for all $f \in \partial H_d$.  If \eqref{ebe-1} fails (i.e. \eqref{ebe-3} holds), then we halt the algorithm (we are in the ``structure'' half of the dichotomy).  Otherwise,
we return to Step 1.
\end{itemize}

Observe that every time we return from Step 2 to Step 1, the quantity
$$ \sum_{e \in H_d} \sum_{E_e \in \B_e} \Energy_{E_e}(\bigvee_{f \in \partial e} \B'_f)$$
increases by at least $\delta^2$.  On the other hand, if this quantity ever increases by more than $|H_d| 2^{2^m} \eps^2 = O_{|J|, m, \eps}(1)$, then
by \eqref{b-card} and the pigeonhole principle \eqref{ebe-1} will necessarily fail.  Since we only
return to Step 1 when \eqref{ebe-1} holds, we see that the algorithm can only iterate at most $O_{|J|, m, \eps, \delta}(1)$ times.  
Thus when we terminate we must have \eqref{ebe-4}.  The claim then folows.
\end{proof}

We now iterate Lemma \ref{dichotomy} to obtain the following preliminary regularity lemma.
Define a \emph{growth function} to be an increasing function $F: \R^+ \to \R^+$ such that $F(x) \geq 1+x$ for all $x$.  

\begin{lemma}[Preliminary regularity lemma]\label{partial-regularity}  Let $V = (J, (V_j)_{j \in J}, d, H_d)$ be a hypergraph system.
For each $e \in H_d$ let $\B_e \subseteq \A_e$ be a $\sigma$-algebra, and suppose that we have the bound
$$ \complexity(\B_e) \leq m \hbox{ for all } e \in H_d$$
for some $m > 0$.  Let $\eps > 0$, and let $F$ be a growth function (possibly depending on $\eps$).  
Then there exists $M > 0$, and for each $f \in \partial H_d$ there exists a pair of $\sigma$-algebras
$\B_f \subseteq \B'_f \subseteq \A_f$ such that we have the estimates
\begin{align}
F(m) \leq M &\leq O_{|J|, \eps, m, F}(1) \label{M-bound} \\
\complexity( \B_f ) &\leq M \hbox{ for all } f \in \partial H_d \label{coarse-complex} \\
\Energy_{E_e}(\bigvee_{f \in \partial e} \B'_f) - \Energy_{E_e}(\bigvee_{f \in \partial e} \B_f) &\leq \eps^2 
\hbox{ for all } e \in H_d, E_e \in \B_e\label{coarse-fine} \\
\Delta_e( E_e | \bigvee_{f \in \partial e} \B'_f ) &\leq \frac{1}{F(M)} \hbox{ for all } e \in H_d, E_e \in \B_e \label{fine-accurate}
\end{align}
\end{lemma}

\begin{remark}
Lemma \ref{partial-regularity} provides a coarse low-order approximation $(\B_f)_{f \in \partial H_d}$ and
a fine low-order approximation $(\B'_f)_{f \in \partial H_d}$ 
to the high-order $\sigma$-algebras $(\B_e)_{e \in H_d}$.  The coarse approximation has bounded complexity,
the fine approximation is close to the coarse approximation in an $L^2$ sense, and the high order $\sigma$-algebras are pseudorandom
with respect to the fine approximation.  The key point here is that the discrepancy control on the fine approximation given by \eqref{fine-accurate} is superior to the complexity control on the coarse approximation given by \eqref{coarse-complex} by an \emph{arbitrary} growth function $F$.  If one were to try to use a single approximation instead of a pair of coarse and fine approximations, it appears impossible to obtain such a crucial gain.
\end{remark}

\begin{proof}  We perform the following iteration.

\begin{itemize}
\item Step 0.  Initialize $\B_f = \{ \emptyset, V_J\}$ to be the trivial $\sigma$-algebra for all $f \in \partial H_d$, thus $\B_f$ has complexity 0
initially.

\item Step 1.  Set $M := \max(F(m), \sup_{f \in \partial H_d} \complexity(\B'_f))$, and $\delta := 1/F(M)$.
We apply Lemma \ref{dichotomy}, and end up in either the randomness or structure half of the dichotomy.  In either
case we generate $\sigma$-algebras $\B_f \subseteq \B'_f \subseteq \A_f$ for each $f \in \partial H_d$.

\item Step 2.  If we are in the randomness half of the dichotomy, we terminate the algorithm.  Otherwise, if we are in the structure half of
the dichotomy, we replace $\B_f$ with $\B'_f$ for each $f \in \partial H_d$, and return to Step 1.
\end{itemize}

Observe that every time we return from Step 2 to Step 1, the quantity
$$ \sum_{e \in H_d} \sum_{E_e \in \B_e} \Energy_{E_e}(\bigvee_{f \in \partial e} \B_f)$$
increases by at least $\eps^2$.  On the other hand, this quantity is non-negative and does not exceed $|H_d| 2^{2^m} = O_{|J|,m}(1)$, thanks to \eqref{b-card}.
Thus this algorithm terminates after $O_{|J|, m, \eps}(1)$ steps.  By \eqref{ebe-4}, we see that at each of these steps, the quantity $M$
increases to be at most $M + O_{J, m, \eps, F(M)}(1)$, while initially $M$ is equal to $F(m)$.  Thus at the end of the algorithm we have 
\eqref{M-bound} as desired.  The remaining claims \eqref{coarse-complex}, \eqref{coarse-fine}, \eqref{fine-accurate} follow from
construction (and \eqref{ebe-1}, \eqref{ebe-2}).
\end{proof}

\begin{remark} Lemma \ref{partial-regularity} already implies the Szemer\'edi regularity lemma in its usual form (and with
the usual tower-exponential bounds); see \cite{tao:regularity} for further discussion.
The above lemma is also similar in spirit to the modern regularity lemmas that appear for instance in \cite{rs} (except for an issue of obtaining regularity at all orders less than $d$, which we shall address in Lemma \ref{full-regularity} below).  
In such lemmas, the objective is not to obtain a partition for which the original graph or hypergraph is regular, but instead to obtain a partition for which a \emph{modified} graph or hypergraph is \emph{very} regular, where the modification consists of adding or subtracting a small number of edges.  The analogue of such a modification in our context is the decomposition
$$ 1_{E_e} = F_{\operatorname{regular}} + F_{\operatorname{small}}$$
where
$$ F_{\operatorname{regular}} := \E( 1_{E_e} | \bigvee_{f \in \partial e} \B_f ) + (1_{E_e} - \E(1_{E_e} | \bigvee_{f \in \partial e} \B'_f))$$
and
$$ F_{\operatorname{small}} := \E(1_{E_e} | \bigvee_{f \in \partial e} \B'_f) - \E(1_{E_e} | \bigvee_{f \in \partial e} \B_f).$$
The function $F_{\operatorname{small}}$ is small thanks to \eqref{coarse-fine} and \eqref{pythagoras}.  Now consider $F_{\operatorname{regular}}$.  On a typical
atom of $\bigvee_{f \in\partial e} \B_f$, the first term is constant, and the second term is going to be very pseudorandom (have small correlation
with sets of the form $\bigcap_{f \in \partial e} E_f$ for $E_f \in \A_f$) thanks to \eqref{fine-accurate} and \eqref{Psie}.
\end{remark}

Lemma \ref{partial-regularity} 
regularizes the $\sigma$-algebras $\B_e$ on the $d$-uniform hypergraph $H_d$ in terms of $\sigma$-algebras $\B_f$, $\B'_f$ on
the $(d-1)$-uniform hypergraph $\partial H_d$.  However it does not regularize the $\sigma$-algebras on $\partial H_d$.  This 
can be accomplished by one final iteration, which gives our final regularity lemma (which is essentially the same lemma\footnote{In contrast, the earlier regularity lemmas of 
Chung \cite{chung} and Frankl-Rodl \cite{frankl} are closer to Lemma \ref{partial-regularity}, with $\partial H_d$ generalized
to $\partial^l H_d$ for any fixed $l$.  The case $l=d-1$ in particular is essentially a routine generalization of the ordinary regularity lemma and appears to have been folklore for quite some time.}  as that in \cite{gowers-hyper}, \cite{rodl}, or \cite{rs}).

\begin{lemma}[Full regularity lemma]\label{full-regularity} Let $V = (J, (V_j)_{j \in J}, d, H_d)$ be a hypergraph system, and define
the $j$-uniform hypergraphs $H_j$ for all $0 \leq j < d$ recursively backwards from $j=d$ by the formula $H_j := \partial H_{j+1}$.
(In particular, if $H_d$ is non-empty then we have $H_0 = \{\emptyset\}$.)  
For all $e \in H_d$ let $\B_e \subseteq \A_e$ be a $\sigma$-algebra, and suppose that we have the bound
$$ \complexity(\B_e) \leq M_d \hbox{ for all } e \in H_d$$
for some $M_d > 0$.  Let $F$ be a growth function.  Then there exists numbers
\begin{equation}\label{growth-cond}
M_d \leq F(M_d) \leq M_{d-1} \leq F(M_{d-1}) \leq \ldots \leq M_0 \leq F(M_0) \leq O_{|J|, M_d, F}(1)
\end{equation}
and for each $0 \leq j < d$ and $f \in H_j$ there exist $\sigma$-algebras
$\B_f \subseteq \B'_f \subseteq \A_f$, such that we have the estimates
\begin{align}
\complexity( \B_f ) &\leq M_j \hbox{ for all } 0 \leq j < d, f \in H_j \label{coarse-complex-2} \\
\Energy_{E_e}(\bigvee_{f \in \partial e} \B'_f) - \Energy_{E_e}(\bigvee_{f \in \partial e} \B_f) &\leq \frac{1}{F(M_j)^2} 
\hbox{ for all } 1 \leq j \leq d, e \in H_j, E_e \in \B_e\label{coarse-fine-2} \\
\Delta_e( E_e | \bigvee_{f \in \partial e} \B'_f ) &\leq \frac{1}{F(M_0)} \hbox{ for all } 1 \leq j \leq d, e \in H_j, E_e \in \B_e.
 \label{fine-accurate-2}
\end{align}
\end{lemma}

\begin{remark}  At every order $0 \leq j \leq d$, Lemma \ref{full-regularity} gives coarse and fine approximations $(\B_f)_{f \in H_{j-1}}$, $(\B'_f)_{f \in H_{j-1}}$
at the $(j-1)$-uniform level to the $\sigma$-algebras $(\B'_e)_{e \in H_j}$ at the $j$-uniform level.  As one goes down in order, the
$\sigma$-algebras rapidly become more complex\footnote{At the zeroth order $j=0$, all $\sigma$-algebras have complexity zero, but this is a degenerate exception to the above general rule.} (though lower order, of course).  However, the bounds in \eqref{coarse-fine-2} and
\eqref{fine-accurate-2} will keep apace with this growth in complexity (see \cite{rs} for some related discussion concerning the desirability
of having the constants grow along such a hierarchy).  Indeed the bound \eqref{fine-accurate-2} is extremely strong,
as $F(M_0)$ dominates all the other quantities which appear in the above lemma; it is effectively as if the fine approximation was perfectly accurate
(so that $1_{E_e}$ is approximable by $\E(1_{E_e} |\bigvee_{f \in \partial e} \B'_f )$ with only negligible error).  The main remaining
difficulty when using this lemma is to exploit the estimate \eqref{coarse-fine-2} measuring the gap between the coarse and fine approximations;
one has to take some care here because 
the error bound $1/F(M_j)^2$ here safely exceeds the complexity\footnote{We will only need to bound the complexity of the coarse algebras $\B_e$.  Some (very weak) bounds on the complexity of the fine algebras $\B'_e$ are available but they seem to be useless for applications and so we have not stated them explicitly here.} of the higher-order objects $(\B_e)_{e \in H_j}$, but not that of
the lower-order objects $(\B_e)_{e \in H_{j-1}}$.  
\end{remark}

\begin{proof}  We induct on $d$ (keeping $J$ fixed); the implicit constants in \eqref{growth-cond} will change when one does this, but the induction
will only run for at most $|J|$ steps and so this will not cause a difficulty.  
When $d=0$ the claim is trivial (and the claim \eqref{coarse-complex-2} has an enormous amount of room available!) 
so assume that
$d \geq 1$ and the claim has already been proven for all smaller $d$.
We will need a growth function $F^{\operatorname{fast}}$ to be chosen later; as the name suggests, this function will grow substantially faster than $F$,
in particular we assume $F^{\operatorname{fast}}(n) \geq F(n)$ for all $n$.
Applying Lemma \ref{partial-regularity}  with $m$ equal to $M_d$, with $\eps$ equal to $1/F(M_d)$, and
the growth function $F^{\operatorname{fast}}$, we can create $\sigma$-algebras $\B_f \subseteq \B'_f \subseteq \A_f$ for all
$f \in H_{d-1}$ and a quantity $M_{d-1}$ such that
\begin{align}
F(M_d) \leq F^{\operatorname{fast}}(M_d) \leq M_{d-1} &\leq O_{|J|, \eps, M_d, F^{\operatorname{fast}}}(1) 
= O_{|J|, M_d, F, F^{\operatorname{fast}}}(1) \label{M-bound-0} \\
\complexity( \B_f ) &\leq M_{d-1} \hbox{ for all } f \in H_{d-1} \nonumber \\
\Energy_{E_e}(\bigvee_{f \in \partial e} \B'_e) - \Energy_{E_e}(\bigvee_{f \in \partial e} \B_f) &\leq \frac{1}{F(M_d)^2} 
\hbox{ for all } e \in H_d, E_e \in \B_e\nonumber \\
\Delta_e( E_e | \bigvee_{f \in \partial e} \B'_f ) &\leq \frac{1}{F^{\operatorname{fast}}(M_{d-1})} \hbox{ for all } e \in H_d, E_e \in \B_e. \label{fine-accurate-0}
\end{align}
Now we apply the induction hypothesis with $d$ replaced by $d-1$, and $H_d$ replaced by $H_{d-1}$.  This generates numbers
\begin{equation}\label{mmm}
M_{d-1} \leq F(M_{d-1}) \leq \ldots \leq M_0 \leq F(M_0) \leq O_{|J|, M_{d-1}, F}(1)
\end{equation}
and for each $0 \leq j < d-1$ and $f \in H_j$ there exist $\sigma$-algebras
$\B_f \subseteq \B'_f \subseteq \A_f$, such that we have the estimates
\begin{align*}
\complexity( \B_f ) &\leq M_j \hbox{ for all } 0 \leq j < d-1, f \in H_j  \\
\Energy_{E_e}(\bigvee_{f \in \partial e} \B'_e) - \Energy_{E_e}(\bigvee_{f \in \partial e} \B_f) &\leq \frac{1}{F(M_j)^2} 
\hbox{ for all } 1 \leq j \leq d-1, e \in H_j, E_e \in \B_e\\
\Delta_e( E_e | \bigvee_{f \in \partial e} \B'_f ) &\leq \frac{1}{F(M_0)} \hbox{ for all } 1 \leq j \leq d-1, e \in H_j, E_e \in \B_e.
\end{align*}
Comparing this with the conclusion of Lemma \ref{full-regularity}, we see that we can obtain all the claims we need except for
\eqref{fine-accurate-2} when $j=d$, as well as the final bound in \eqref{growth-cond}.  To obtain \eqref{fine-accurate-2}, we see from
\eqref{fine-accurate-0} that it would suffice to ensure that
$$ F^{\operatorname{fast}}(M_{d-1}) \geq F(M_0).$$
But since $F(M_0) = O_{|J|,M_{d-1}, F}(1)$, this can be achieved simply by choosing the growth function $F^{\operatorname{fast}}$ to be sufficiently large and
rapidly increasing depending on $F$ and $|J|$.  By \eqref{mmm}, \eqref{M-bound-0}, we then have
$$ F(M_0) = O_{|J|,M_{d-1}, F}(1) = O_{|J|, M_d, F, F^{\operatorname{fast}}}(1) = O_{|J|, M_d, F}(1)$$
and the claim \eqref{growth-cond} follows.
\end{proof}

\begin{remark} The dependence of constants here is quite terrible.  Typically $F$ will be an exponential function.  In the graph case
$d=2$ one can take $M_0$ to be a tower of exponentials, whose height is bounded by some polynomial of $F(M_2)$; a 
modification of the arguments in \cite{gowers-sz}
shows that this tower bound is essentially best possible.  However, for $d=3$, both $M_0$ and $M_1$ will be an \emph{iterated} tower of exponentials of iterated height equal to a polynomial in $F(M_3)$, basically because of the need for $F^{\operatorname{fast}}$ to exceed the bounds one obtains from the $d=2$ case.  The situation of course gets even worse for larger values of $d$, though 
for any fixed $d$ the bounds are still primitive recursive.  As stated earlier, the complexity bounds for the fine approximations
$\B'_f$ will be even worse than this, perhaps by yet another layer of iteration.
Nevertheless, this regularity lemma is still sufficient for applications in which one is willing to have qualititative control only on the error
terms (e.g. $o(1)$ type bounds) rather than quantitative control.  (As we shall see in \cite{tao-multiprime}, obtaining infinitely many constellations in the Gaussian primes will be one such application.)  In view of recent results on effective bounds on Szemer\'edi-type theorems (see e.g. \cite{gowers}, \cite{shkredov}) it seems quite possible that these very rapid bounds, while perhaps necessary in order to have a regularity lemma, are not needed for the hypergraph removal lemma.
\end{remark}

\section{Statement of counting lemma}

As is customary in these arguments, the regularity lemma must be complemented with a counting lemma in order for it to be applicable
to proving results such as Theorem \ref{main-2}.  In the $\sigma$-algebra language, the setup is as follows.  Suppose
we start with $\sigma$-algebras $(\B_e)_{e \in H_d}$ as in the hypotheses of Lemma \ref{full-regularity}.
Then, among other things, this lemma yields further $\sigma$-algebras $(\B_e)_{e \in H_j}$ for $0 \leq j < d$, each of which has some complexity bound.
Combining all of these $\sigma$-algebras together, one obtains a somewhat large (but still bounded complexity) $\sigma$-algebra
$ \bigvee_{e \in H} \B_e$, where $H := \bigcup_{0 \leq j \leq d} H_j$.  In particular, if $E_e$ are sets in $\B_e$ for all 
$e \in H_d$, then $\bigcap_{e \in H_d} E_e$ is the union of atoms in $\bigvee_{e \in H} \B_e$.  Here, of course, an atom of a $\sigma$-algebra $\B$
is a non-empty set in $\B$ of minimal size; since the ambient space $V_J$ is finite, every point is contained in exactly one atom of $\B$.

Roughly speaking, the counting lemma we give below (Lemma \ref{count-lemma})
gives a formula for computing the probability of atoms in $\bigvee_{e \in H} \B_e$, or at least those atoms which are ``good''.
It can be informally described as follows.  For each $e \in H$, let $A_e$ be an atom of $\B_e$, thus
$\bigcap_{e \in H} A_e$ will be an atom of $\bigvee_{e \in H} \B$ (if it is non-empty).  The counting lemma then says that under most circumstances we have the approximate formula\footnote{The reader may wish to interpret $\E(1_A)$ as being the ``probability'' of the ``event'' $A$, thus
for instance $\E( \prod_{e \in H} 1_{A_e})$ is the probability of the joint event $\bigcap_{e \in H} A_e$.  Similarly, many of the arguments in the
sequel also have a strongly probabilistic flavour.}
\begin{equation}\label{counting}
\E( \prod_{e \in H} 1_{A_e} ) \approx \prod_{e \in H} \E( 1_{A_e} | \bigcap_{f \in \partial e} A_f )
\end{equation}
where we use $\E(f|A)$ to denote the conditional expectation
$$ \E(f|A) := \frac{1}{|A|} \sum_{x \in A} f(x).$$
This can be viewed as an assertion that higher order atoms $A_e$ are approximately independent of each other, conditioning on lower order
atoms $A_f$, although a precise formulation of this heuristic is somewhat difficult to quantify.  In particular, if we remove those ``bad''
atoms $\bigcap_{e \in H} A_e$
for which $\E( 1_{A_e} | \bigcap_{f \in \partial e} A_f )$ is small for at least one $e \in H$, then all the remaining non-empty atoms will have fairly large size.  Thus if the set $\bigcap_{e \in H} E_e$ has very small size, then after removing all the bad atoms we expect
this set to in fact be empty.  This is the strategy behind proving Theorem \ref{main-2}.

We now formalize the above discussion.  We begin by describing the good atoms.  Informally speaking, the good atoms are going to be those which are fairly large (at all orders) and also fairly regular (at all orders).  This is consistent with previous 
experience with counting lemmas (say in the graph case), in which one must first throw away all cells of the partition which are 
too small (or have too few edges), as well as all pairs of cells for which the graph is irregular, before one can obtain a useful estimate for (say) the number of triangles in a graph.

\begin{definition}[Good atoms]\label{good-def}  Let the notation, assumptions, and conclusions be as in Lemma \ref{full-regularity}, and let $H := \bigcup_{0 \leq j \leq d} H_j$.  Let $\bigcap_{e \in H} A_e$ be a (possibly empty) atom of $\bigvee_{e \in H} \B_e$, where for each $e \in H$, $A_e$ is an atom of $\B_e$.  We say that this atom is \emph{good} if 
for all $0 \leq j \leq d$ and $e \in H_j$ we have the largeness estimates
\begin{equation}\label{e-large}
\E( 1_{A_e} \prod_{f \in \partial e} 1_{A_f} ) \geq \frac{1}{\log F(M_j)} \E(\prod_{f \in \partial e} 1_{A_f})
\end{equation}
as well as the regularity estimates
\begin{equation}\label{e-regularity}
\E\left( \bigl|\E( 1_{A_e} | \bigvee_{f \in \partial e} \B'_f ) - \E( 1_{A_e} | \bigvee_{f \in \partial e} \B_f )\bigr|^2 \prod_{f \subsetneq e} 1_{A_f} \right)
\leq \frac{1}{F(M_j)} \E( \prod_{f \subsetneq e} 1_{A_f} ).
\end{equation}
\end{definition}

\begin{remark}
While the definition of a good atom allows for $\bigcap_{e\in H} A_e$ to be empty, the counting lemma we prove below will show that in fact good atoms are always non-empty (assuming $F$ is sufficiently rapid). The reader should not take the logarithmic factor in \eqref{e-large} too seriously; the point is that $\log F(M_j)$ is smaller than any power of $F(M_j)$ but still much larger
than any given function of $M_j$.
\end{remark}

One can easily verify that most atoms are good in the following sense.  For any $0 \leq j \leq d$, $e \in H_j$, and any atom 
$A_e$ of $\B_e$, let $B_{e,A_e}$ be the union of all the sets $\bigcap_{f \subsetneq e} A_f$ for which \eqref{e-large} or \eqref{e-regularity} fails.  We remark for future reference that the set $B_{e,A_e}$ lies in $\bigvee_{f \subsetneq e} \B_f$.  
Note also that if the atom $\bigcap_{e \in H} A_e$ is not good, then there exists $e \in H$ such that
$\bigcap_{e' \in H} A_{e'} \subseteq A_e \cap B_{e,A_e}$.

\begin{lemma}[Most atoms are good]\label{good-lots}  Let the notation, assumptions, and conclusions be as in Lemma \ref{full-regularity}
and Definition \ref{good-def}.
For any $0 \leq j \leq d$, $e \in H_j$, and any atom $A_e$ of $\B_e$,
we have $\E(1_{A_e} 1_{B_{e,A_e}}) = O(1 / \log F(M_j))$.
\end{lemma}

\begin{proof}  Consider the contribution to $\E( 1_{A_e} 1_{B_{e,A_e}} )$ from the case where \eqref{e-large} fails.  This contribution is bounded by\footnote{Note that \eqref{e-large} depends only on those $A_f$ for which $f \in \partial e$, as opposed to the larger class of events $A_f$ for which $f \subsetneq e$.}
$$ \sum_{(A_f)_{f \in \partial e} \hbox{\scriptsize atoms in } (\B_f)_{\partial e}: \hbox{\scriptsize \eqref{e-large} fails}}
\E( 1_{A_e} \prod_{f \in \partial e} 1_{A_f} )$$
which by failure of \eqref{e-large} is bounded by
$$ \leq \sum_{(A_f)_{f \in \partial e} \hbox{\scriptsize atoms in } (\B_f)_{\partial e}} \frac{1}{\log F(M_j)} \E( \prod_{f \in \partial e} 1_{A_f} )
= \frac{1}{\log F(M_j)}.$$
Next, consider the contribution to $\E( 1_{A_e} 1_{B_{e,A_e}})$ arising from the case when \eqref{e-regularity} fails.  The total contribution
of this case is 
$$ \sum_{(A_f)_{f \subsetneq e}: \hbox{\scriptsize \eqref{e-regularity} fails}} \E( \prod_{f \subsetneq e} 1_{A_{f}} )$$
which by failure of \eqref{e-regularity} is at most
$$ F(M_j) \sum_{(A_{f})_{f \subsetneq e}} \E\left( \bigl|\E( 1_{A_e} | \bigvee_{f \in \partial e} \B'_f ) - \E( 1_{A_e} | \bigvee_{f \in \partial e} \B_f )\bigr|^2 \prod_{f \subsetneq e} 1_{A_{f}} \right)$$
which in turn is at most
$$ F(M_j) \E\left( |\E( 1_{A_e} | \bigvee_{f \in \partial e} \B'_f ) - \E( 1_{A_e} | \bigvee_{f \in \partial e} \B_f )|^2 \right).$$
But by \eqref{pythagoras}, \eqref{coarse-fine-2} we have
$$ \E\left( |\E( 1_{A_e} | \bigvee_{f \in \partial e} \B'_f ) - \E( 1_{A_e} | \bigvee_{f \in \partial e} \B_f )|^2 \right) \leq \frac{1}{F(M_j)^2}.$$
Combining all of these estimates, the claim follows.
\end{proof}

We can now state the counting lemma; closely related results appear in the work of Gowers \cite{gowers}, Nagle, R\"odl, and Schacht \cite{nrs},
and R\"odl and Schacht \cite{rs}.

\begin{lemma}[Counting lemma]\label{count-lemma}  Let the notation, assumptions, and conclusions be as in Lemma \ref{full-regularity} and
Definition \ref{good-def}, and let $H := \bigcup_{0 \leq j \leq d} H_j$.  Let $\bigcap_{e \in H} A_e$ be a good atom of $\bigvee_{e \in H} \B_e$.  Then, if the growth function $F$ is sufficiently rapid depending on $|J|$, we have that $\bigcap_{e \in H} A_e$ is non-empty, and more precisely
$$\E( \prod_{e \in H} 1_{A_e} ) = (1 + o_{M_d \to \infty; |J|}(1)) \prod_{e \in H} \E( 1_{A_e} | \bigcap_{f \in \partial e} A_f )
+ O_{|J|, M_0}\left(\frac{1}{F(M_0)}\right)
$$
(compare with \eqref{counting}).  
\end{lemma}

This lemma is a little lengthy (though straightforward) to prove, and we defer it to the next section.  Let us assume it for now, and conclude
the proof of Theorem \ref{main-2}.

\begin{proof}[of Theorem \ref{main-2} assuming Lemma \ref{count-lemma}]  Let $V = (J, (V_j)_{j \in J}, d, H_d)$, 
$(E_e)_{e \in H_d}$, $\delta$ be as in Theorem
\ref{main-2}.  We define $H_j$ recursively for $0 \leq j < d$ by setting $H_j := \partial H_{j+1}$, and then set $H := \bigcup_{0 \leq j \leq d} H_j$.
For any $e \in H_d$ we set $\B_e := \B(E_e)$,
thus each $\B_e$ has complexity at most 1.  Let $M_d \geq 1$ be a quantity to be chosen later, and let $F$ be a growth function depending on $|J|$
(but not on $\delta$) to be chosen later.  We apply the regularity lemma, Lemma \ref{full-regularity}, to obtain quantities \eqref{growth-cond}
and $\sigma$-algebras $\B_f \subseteq \B'_f \subseteq \A_f$ for all $f \in H$.

Suppose that $\bigcap_{e \in H} A_e$ is a (possibly empty) atom of $\bigvee_{e \in H} \B_e$ such that $A_e = E_e$ for $e \in H_d$.  If this atom is good, then by the counting Lemma (Lemma \ref{count-lemma}) and Definition \ref{good-def}
we have
$$
\E( 1_{\bigcap_{e \in H} A_e} ) = (1 + o_{M_d \to \infty; |J|}(1)) \prod_{0 \leq j \leq d} \prod_{e \in H_j} \frac{1}{F(M_j)^{1/10}} 
+ O_{|J|, M_0}\left(\frac{1}{F(M_0)}\right),
$$
if $F$ is sufficiently rapid depending on $|J|$.  Using \eqref{growth-cond}, we thus see
that (if $M_d$ is sufficiently large depending on $J$)
$$ \E( 1_{\bigcap_{e \in H} A_e} ) \geq c(|J|, M_d, F)$$
for some $c(|J|, M_d, F) > 0$.  On the other hand, $\bigcap_{e \in H} A_e$ is contained in $\bigcap_{e \in H_d} E_e$, which has density
at most $\delta$ by the hypothesis \eqref{E-dens}.  Thus if $\delta$ is sufficiently small depending on $|J|$, $M_d$,
$F$, we see that no atom $\bigcap_{e \in H} A_e$ with $A_e = E_e$ for $e \in H_d$ can possibly be good.

Now let $B_{e,A_e}$ be as in Lemma \ref{good-lots}.  Let us define
$$ E'_e := V_J \backslash \bigl( B_{e, E_e} \cup \bigcup_{f \subsetneq e} \bigcup_{A_f} A_f \cap B_{f, A_f} \bigr)$$
for all $e \in H_d$, where for brevity we adopt the convention that $A_f$ is always understood to range over the atoms of $\B_f$.  Then we observe that $E'_e \in \bigvee_{f \subsetneq e} \B_f$.  The claims \eqref{E-complex}, \eqref{E-meas}
then follow from \eqref{coarse-complex-2}.
Also, from Lemma \ref{good-lots}, \eqref{coarse-complex-2} we see that for any $e \in H_d$,
\begin{align*}
\E( 1_{E_e \backslash E'_e} ) &\leq \E( 1_{E_e} 1_{B_{e,E_e}} ) + 
\sum_{f \subsetneq e} \sum_{A_f} \E( 1_{A_f} 1_{B_{f, A_f}} ) \\
&\leq O(F(M_d)^{-1/10}) + \sum_{0 \leq j < d} \sum_{f \in H_j} \sum_{A_f}  O( 1 / \log F(M_j) ) \\
&\leq O(F(M_d)^{-1/10}) + \sum_{0 \leq j < d} \sum_{f \in H_j} O_{M_j}( 1 / \log F(M_j) ) \\
&\leq \sup_{0 \leq j \leq d} O_{M_j, |J|}(1 / \log F(M_j)).
\end{align*}
If one chooses $F$ sufficiently rapidly growing (depending only on $|J|$), we conclude from \eqref{growth-cond} that we have
$$ \E(1_{E_e \backslash E'_e}) = o_{M_d \to 0; |J|}(1).$$
By choosing $M_d$ sufficiently large depending on $|J|$, and then letting $\delta$ be sufficiently small depending on $M_d$ and $|J|$, we
conclude \eqref{E-error}.

The final thing to verify is \eqref{E-cap}.  To see this, first observe that this set lies in $\bigvee_{f \in H \backslash H_d} \B_f$ and thus is the union of atoms of the form $\bigcap_{f \in H \backslash H_d} A_f$.  Suppose for contradiction that $\bigcap_{e \in H_d} E'_e$
contains a non-empty atom of the form $\bigcap_{f \in H \backslash H_d} A_f$. Set $A_e := E_e$ for $e \in H_d$.  By the preceding discussion
we know that $\bigcap_{e \in H} A_e$ cannot be good, thus there exists an $f' \in H$ such 
that $\bigcap_{g \subsetneq f'} A_g$ lies in $B_{f',A_{f'}}$.
From construction of $H$, there exists $e \in H_d$ which contains $f'$.
But then by definition of $E'_e$, $\bigcap_{f \in H \backslash H_d} A_f$ cannot lie in $E'_e$, contradiction.  Thus 
$\bigcap_{e \in H_d} E'_e$ is empty, which is \eqref{E-cap}, and Theorem \ref{main-2} follows.
\end{proof}

It remains to prove the counting lemma.  This will be accomplished in the next section.

\section{Proof of counting lemma}

We now prove Lemma \ref{count-lemma}.  Fix a good collection $(A_e)_{e \in H}$ of atoms.  We introduce the numbers $p_e \in \R$,
the functions $b_e, c_e: V_J \to \R$, and the sets $A_{<e} \subseteq V_J$ for all $e \in H$ by the formulae	
\begin{align*}
p_e &:= \E( 1_{A_e} | \bigcap_{f \in \partial e} A_f ) \\
b_e &:= \E( 1_{A_e} | \bigvee_{f \in \partial e} \B'_f ) - \E( 1_{A_e} | \bigvee_{f \in \partial e} \B_f ) \\
c_e &:= 1_{A_e} - \E( 1_{A_e} | \bigvee_{f \in \partial e} \B'_f ) \\
A_{<e} &:= \bigcap_{f \subsetneq e} A_f.
\end{align*}
Note that we have not yet shown that $\bigcap_{f \in \partial e} A_f$ is non-empty; for now, let us just assign an arbitrary value 
to $p_e$ (e.g. $p_e = 1$) when $\bigcap_{f \in \partial e} A_f$ is empty.  We thus have the decomposition
\begin{equation}\label{e-decomp}
1_{A_e} = p_e + b_e + c_e 
\end{equation}
on the set $\bigcap_{f \in \partial_e} A_f$.  One should think of the constant $p_e$ as the main term, and the other two
terms as error terms.  The $c_e$ error term will be very easy to handle, whereas the $b_e$ error term will cause
somewhat more difficulty.  Since $(A_e)_{e \in H}$ is good, we have the estimates
\begin{equation}\label{pe-big}
p_e \geq 1 / \log F(M_j) \hbox{ for all } 0 \leq j \leq d \hbox{ and } e \in H_j
\end{equation}
and
\begin{equation}\label{ge-small}
\E( |b_e|^2 1_{A_{<e}} ) \leq F(M_j)^{-1} \E( 1_{A_{<e}} ) \hbox{ for all } 0 \leq j \leq d \hbox{ and } e \in H_j.
\end{equation}
From \eqref{fine-accurate-2} and \eqref{Psie}, we also have
\begin{equation}\label{he-small}
|\E( c_e \prod_{f \in \partial e} 1_{E_f} )| \leq \frac{1}{F(M_0)} \hbox{ whenever } E_f \in \A_f \hbox{ for } f \in \partial e.
\end{equation}

Our objective is to use the above estimates \eqref{e-decomp}, \eqref{pe-big}, \eqref{ge-small}, \eqref{he-small} to conclude that
\begin{equation}\label{ae}
 \E( \prod_{e \in H} 1_{A_e} ) = (1 + o_{M_d \to \infty; |J|}(1)) \prod_{e \in H} p_e + O_{|J|, M_0}(\frac{1}{F(M_0)}).
 \end{equation}
This will be achieved by several applications of the Cauchy-Schwarz and triangle inequalities.  However, there is a certain amount of notational
burden in order to keep track of the expressions in the succesive applications of these inequalities.    It will be convenient to
return to the original sets $(V_j)_{j \in J}$.  We can identify $A_e \in \B_e$ as a subset $\overline{A_e}$ of $V_e = \prod_{j \in e} V_j$, and similarly we can view
the $\A_e$-measurable $b_e$ and $c_e$ as functions $\overline{b_e}$ and $\overline{c_e}$ 
on $V_e$. One can then write \eqref{ae} in the form
\begin{equation}\label{vj-form}
\begin{split}
\frac{1}{\prod_{j \in J} |V_j|} &\sum_{(v_j)_{j \in J} \in \prod_{j \in J} V_j}\ \prod_{e \in H} 1_{\overline{A_e}}\bigl( (v_j)_{j \in e} \bigr) 
\\
&=
\bigl(1 + o_{M_d \to \infty; |J|}(1)\bigr) \prod_{e \in H} p_e + O_{|J|, M_0}\left(\frac{1}{F(M_0)}\right).
\end{split}
\end{equation}
For inductive purposes we will need to generalize\footnote{The basic problem is that we need the Cauchy-Schwarz inequality to eliminate each of the $\overline{b_e}$ factors in turn (using \eqref{ge-small}), but each time we apply this inequality we essentially double the number of free variables that one has to sum or average over.  In particular, one ends up sampling more than one point from 
each vertex class $V_j$, which forces us to leave the probabilistic framework that has been so convenient for us in preceding sections and return to a combinatorial framework.  One could stay in the probabilistic framework using the machinery of tensor products (and conditional tensor products) of probability spaces, but this would introduce even more excessive notation into an already notation-heavy argument and would probably not be helpful to the reader.} 
 this formula.

\begin{definition}[Hypergraph bundle]  A \emph{hypergraph bundle} over $H$ is a hypergraph $G \subseteq 2^K$ on a finite set $K$, together
with a map $\pi: K \to J$ (which we call the \emph{projection map} of the bundle), which is a hypergraph homomorphism (i.e. for each edge $g \in G$, the function $\pi$ is injective on $g$ and $\pi(g) \in H$).  For any $g \subseteq K$, we write $V_g$ for the product set $V_g := \prod_{k \in g} V_{\pi(k)}$.  We say that the bundle is \emph{closed under set inclusion} if whenever $g \in G$ and $g' \subset g$, we have $g' \in G$.
\end{definition}

\begin{remark} From a probabilistic viewpoint, the probability space $V_J$ corresponds to sampling one vertex independently from each of the vertex classes $V_j$ of $V_J$, whereas the more general spaces $V_g$ correspond to the possibility of sampling more than one vertex independently from each of the vertex classes.
\end{remark}

The generalization of the formula \eqref{vj-form} is then

\begin{lemma}[Generalized counting lemma]\label{gencount}  Let $G \subseteq 2^K$ be a hypergraph bundle over $H$ which is closed under set inclusion, 
with projection map $\pi: K \to J$.  Let
$d' := \sup_{g \in G} |g|$ be the order of $G$.  Then, if $F$ is sufficiently rapidly growing depending on $d'$, $|J|$ and $|K|$, we have
\begin{equation}\label{vk-count}
\begin{split}
&\frac{1}{|V_K|} \sum_{(v_k)_{k \in K} \in V_K}\ \prod_{g \in G} 1_{\overline{A_{\pi(g)}}}( (v_k)_{k \in g} ) \\
&=
\bigl(1 + o_{M_d \to \infty; d', |J|, |K|}(1)\bigr) \prod_{g \in G} p_{\pi(g)} + O_{d', |J|, |K|, M_0}\left(\frac{1}{F(M_0)}\right).
\end{split}
\end{equation}
\end{lemma}

Observe that \eqref{vj-form} is the special case of this lemma with $G = H$ (and $K = J$, and $\pi$ being the identity map); note from construction of $H$ that $H$ is automatically closed under set inclusion.

\begin{proof}  We shall use a double induction.  Firstly, we shall induct on the order $d'$ of the bundle $G$.  When $d' = 0$ the claim is
vacuously true (the left-hand side and the main term of the right-hand side is equal to 1), so we may assume $d' \geq 1$ and the claim has already been proven for $d'-1$ and for all choices of hypergraph bundle $G \subseteq 2^K$ which are closed under set inclusion.

Next, we fix $K$ and induct on the quantity $r := |\{ g \in G: |g| = d' \}|$, which is a positive integer between $1$ and $2^{|K|}$.  We thus assume
that the claim has already been proven for all smaller values of $r$ (note that for $r=0$ this follows from the previous induction hypothesis).
The constants may change as we progress in this induction, but since the number of steps in the induction cannot exceed $2^{|K|}$, this will not be a concern.

Let $g_0 \in G$ be such that $|g_0| = d'$.  We use \eqref{e-decomp} to split
\begin{align*}&\prod_{g \in G} 1_{\overline{A_{\pi(g)}}}( (v_k)_{k \in g} ) =\\
&\left[\prod_{g \in G \backslash \{g_0\}} 1_{\overline{A_{\pi(g)}}}( (v_k)_{k \in g} )\right] 
\left( p_{\pi(g_0)} + \overline{b_{\pi(g_0)}}\bigl((v_k)_{k \in g_0}\bigr) + \overline{c_{\pi(g_0)}}\bigl((v_k)_{k \in g_0}\bigr) \right)
\end{align*}
and consider the contribution of the three terms separately.

We first consider the contribution of the $p_{\pi(g)}$ term, which is the main term.  Applying the second induction hypothesis to $G \backslash \{g_0\}$ we see from
\eqref{vk-count} that
\begin{align*}
&\frac{1}{|V_K|} \sum_{(v_k)_{k \in K} \in V_K} \prod_{g \in G \backslash \{g_0\}} 1_{\overline{A_{\pi(g)}}}\bigl( (v_k)_{k \in g} \bigr)
\\
&=
\bigl(1 + o_{M_d \to \infty; d', |J|, |K|}(1)\bigr) \prod_{g \in G \backslash \{g_0\}} p_{\pi(g)} + O_{d', |J|, |K|, M_0}\left(\frac{1}{F(M_0)}\right).
\end{align*}
Multiplying this by the quantity $p_{\pi(g_0)}$, which is between 0 and 1, we see that the contribution of this term to \eqref{vk-count} is
\begin{equation}\label{contrib-1}
(1 + o_{M_d \to \infty; d', |J|, |K|}(1)) \prod_{g \in G} p_{\pi(g)} + O_{d', |J|, |K|, M_0}(\frac{1}{F(M_0)}).
\end{equation}

Next we consider the $\overline{c_{\pi(g_0)}}$ term.  We split $V_K = V_{g_0} \times V_{K \backslash g_0}$.  Let us temporarily freeze the values of $v_k$
for $k \in K \backslash g_0$, and consider the expression
$$
\frac{1}{|V_{g_0}|} \sum_{(v_k)_{k \in g_0} \in V_{g_0}} 
\left[ \prod_{g \in G \backslash \{g_0\}} 1_{\overline{A_{\pi(g)}}}\bigl( (v_k)_{k \in g} \bigr) \right]
\overline{c_{\pi(g_0)}}\bigl((v_k)_{k \in g_0} \bigr).$$
Observe that for each $g \in G \backslash \{g_0\}$, we have $g \neq g_0$ and $|g| \leq d' = |g_0|$.  Thus $g \cap g_0$ is a proper subset of $g_0$, and
thus there exists an element of $\partial g_0$ which contains $g \cap g_0$.  Thus one can rewrite the product $\prod_{g \in G \backslash \{g_0\}} 1_{\overline{A_{\pi(g)}}}\bigl( (v_k)_{k \in g} \bigr)$ in the form
$$ \prod_{f \in \partial g_0} 1_{E_f}\bigl( (v_k)_{k \in \pi(f)} \bigr)$$
for some sets $E_f \subseteq V_f$ whose exact form is not important here (we allow the $E_f$ to depend on the frozen $v_k$).  Applying \eqref{he-small}, we conclude that
$$
\left|\frac{1}{|V_{g_0}|} \sum_{(v_k)_{k \in g_0} \in V_{g_0}} 
\left[ \prod_{g \in G \backslash \{g_0\}} 1_{\overline{A_{\pi(g)}}}\bigl( (v_k)_{k \in g} \bigr) \right]
\overline{c_{\pi(g_0)}}\bigl((v_k)_{k \in g_0} \bigr)\right| \leq 1/F(M_0).$$
Averaging this over all choices of the frozen variables $k \in K \backslash g_0$, we conclude that the contribution of this term to \eqref{vk-count}
is at most 
\begin{equation}\label{fm0}
1/F(M_0).
\end{equation}

Finally we consider the contribution of the $\overline{b_{\pi(g_0)}}$ term, which is the most difficult from a notational viewpoint to handle, 
mainly because of the need to invoke the Cauchy-Schwarz inequality.  We expand this contribution as
$$
\frac{1}{|V_K|} \sum_{(v_k)_{k \in K} \in V_K} 
\left[ \prod_{g \in G \backslash \{g_0\}} 1_{\overline{A_{\pi(g)}}}\bigl( (v_k)_{k \in g} \bigr) \right]
 \overline{b_{\pi(g_0)}}\bigl((v_k)_{k \in g_0}\bigr).$$
We take absolute values and discard\footnote{This discarding step is important as it lowers the total order of the expression being computed, which compensates for a certain doubling of the hypergraph bundle which shall occur shortly when we apply Cauchy-Schwarz.  We can get away with this step because the smallness of $b_{\pi(g_0)}$, as given by \eqref{ge-small}, safely dominates any loss we absorb by discarding these high-order factors.} the bounded factors $1_{\overline{A_{\pi(g)}}}( (v_k)_{k \in g} )$ with $|g| = d'$, to estimate this expression by
$$
O\left( \frac{1}{|V_K|} \sum_{(v_k)_{k \in K} \in V_K} \left[ \prod_{g \in G_{\subsetneq g_0} \cup G'} 1_{\overline{A_{\pi(g)}}}\bigl( (v_k)_{k \in g} \bigr) \right]
 \bigl|\overline{b_{\pi(g_0)}}\bigl((v_k)_{k \in g_0}\bigr)\bigr| \right)$$
 where $G_{\subsetneq g_0} := \{ g: g \subsetneq g_0 \}$ and $G' := \{g \in G \backslash G_{\subsetneq g_0}: |g| \leq d'-1 \}$.  We factorize this
 as
\begin{equation}\label{precauchy}
\begin{split}
O\biggl( \frac{1}{|V_{g_0}|} \sum_{(v_k)_{k \in g_0} \in V_{g_0}}& 
\left[ \prod_{g \in G_{\subsetneq g_0}} 1_{\overline{A_{\pi(g)}}}\bigl( (v_k)_{k \in g} \bigr) \right] \left|\overline{b_{\pi(g_0)}}\bigl((v_k)_{k \in g_0}\bigr)\right|\\
& 
\left[ \frac{1}{|V_{K \backslash g_0}|} \sum_{(v_k)_{k \in K \backslash g_0} \in V_{K \backslash g_0}}
\prod_{g \in G'} 1_{\overline{A_{\pi(g)}}}\bigl( (v_k)_{k \in g} \bigr) \right] \biggr).
\end{split}
\end{equation}

On the other hand, from \eqref{ge-small} we have
\begin{align*}
\frac{1}{|V_{g_0}|} \sum_{(v_k)_{k \in g_0} \in V_{g_0}} &
\left[\prod_{g \in G_{\subsetneq g_0}} 1_{\overline{A_{\pi(g)}}}\bigl( (v_k)_{k \in g} \bigr)\right]\\
& \bigl|\overline{b_{\pi(g_0)}}\bigl((v_k)_{k \in g_0}\bigr)\bigr|^2
\leq \frac{1}{F(M_{d'})} \E(1_{\overline{A_{<\pi(g_0)}}}),
\end{align*}
and hence by Cauchy-Schwarz we can estimate \eqref{precauchy} by

\begin{equation}\label{post-cauchy}
\begin{split}
O\Biggl( F(M_{d'})^{-1/2} &\E(1_{\overline{A_{<\pi(g_0)}}})^{1/2}
\biggl(\frac{1}{|V_{g_0}|} \sum_{(v_k)_{k \in g_0} \in V_{g_0}} 
\left[\prod_{g \in G_{\subsetneq g_0}} 1_{\overline{A_{\pi(g)}}}\bigl( (v_k)_{k \in g} \bigr)\right]\\  
&\left[ \frac{1}{|V_{K \backslash g_0}|} \sum_{(v_k)_{k \in K \backslash g_0} \in V_{K \backslash g_0}}
\prod_{g \in G'} 1_{\overline{A_{\pi(g)}}}\bigl( (v_k)_{k \in g} ) \bigr) \right]^2 \biggr)^{1/2} \Biggr.
\end{split}
\end{equation}
From the first induction hypothesis we have
$$ \E(1_{\overline{A_{<\pi(g_0)}}}) = \bigl(1 + o_{M_d \to \infty; d', |J|}(1)\bigr) \prod_{g \in G_{\subsetneq g_0}} p_{\pi(g)} + O_{d', |J|, M_0}\left(\frac{1}{F(M_0)}\right)$$
and thus
\begin{equation}\label{pag}
\E(1_{\overline{A_{<\pi(g_0)}}}) = O_{M_d, d', |J|}(\prod_{g \in G_{\subsetneq g_0}} p_{\pi(g)}) + O_{d', |J|, M_0}\left(\frac{1}{F(M_0)}\right).
\end{equation}
Now we estimate the expression in parentheses in \eqref{post-cauchy}.  As we shall see, this expression can be rewritten in a form
which can be handled by the induction hypothesis, but with the hypergraph bundle $G$ replaced by a hypergraph of approximately
twice the size (roughly speaking, we throw away all edges of top order $d'$, and double all the remaining edges that are not
contained in $G_{\subsetneq g_0}$).  It is this doubling which forces us to work with a generalized counting lemma\footnote{There is a possible alternate approach which avoids the Cauchy-Schwarz inequality, and hence the need to work with hypergraph bundles.  One can attempt to use the lower-order induction hypothesis to show some uniform distribution properties concerning the intersections of the lower-order atoms with each other, in order that the contribution of the $b_{g_0}$ error be shown to be negligible.  A model example of such a statement, in the graph setting, would be the assertion that in an $\eps$-regular graph $H$, the number of copies of a fixed small graph $G$ in $H$, with one edge specified to be $(x,y)$, is usually close to a fixed quantity independent of $x$ and $y$, except for a small number of exceptional pairs $(x,y)$.  We will not pursue such an alternate approach here.}
 rather than the original counting lemma.

Let $\tilde K = K \oplus_{g_0} K$ be the set $K \times \{0,1\}$, with the 
elements $(k,0)$ and $(k,1)$ identified for all $k \in g_0$.  There is an obvious projection $\phi: \tilde K \mapsto K$, and hence a map $\pi \circ \phi: \tilde K \to H$.  On $\tilde K$ we also place a hypergraph bundle $\tilde G$, defined as the set $\{ g \times \{i\}: g \in G_{\subsetneq g_0} \cup G', i \in 1,2\}$; note that $g \times \{0\}$ and $g \times \{1\}$ will be identified when $g \in G_{\subsetneq g_0}$.  From the definitions we observe that
\begin{align*}
&\frac{1}{|V_{g_0}|} \sum_{(v_k)_{k \in g_0} \in V_{g_0}} 
\left[\prod_{g \in G_{\subsetneq g_0}} 1_{\overline{A_{\pi(g)}}}\bigl( (v_k)_{k \in g} \bigr) \right]
\left[ \frac{1}{|V_{K \backslash g_0}|} \sum_{(v_k)_{k \in K \backslash g_0} \in V_{K \backslash g_0}}
\prod_{g \in G'} 1_{\overline{A_{\pi(g)}}}\bigl( (v_k)_{k \in g} \bigr) \right]^2  \\
&= \frac{1}{|V_{\tilde K}|} \sum_{(v_{\tilde k})_{\tilde k \in \tilde K} \in V_{\tilde K}}
\prod_{\tilde g \in \tilde G} 1_{\overline{A_{\pi \circ \phi(\tilde g)}}} \bigl( (v_{\tilde k})_{\tilde k \in \tilde g} \bigr).
\end{align*}
Applying the first induction hypothesis, we can write this expression as
\begin{equation}\label{moo}
 (1 + o_{M_d \to \infty; d', |J|, |K|}(1)) \prod_{\tilde g \in \tilde G} p_{\pi \circ \phi(\tilde g)} + O_{d', |J|, |K|, M_0}\left(\frac{1}{F(M_0)}\right).
 \end{equation}
By the definition of $\tilde G$, we can write
$$ \prod_{\tilde g \in \tilde G} p_{\pi \circ \phi(\tilde g)}  = \prod_{g \in G_{\subsetneq g_0}} p_{\pi(g)} \times 
[\prod_{g \in G'} p_{\pi(g)}]^2$$
and thus by \eqref{pe-big} and \eqref{growth-cond} we can rewrite \eqref{moo} as
$$ O_{M_d, d', |J|, |K|}\left(\prod_{g \in G_{\subsetneq g_0}} p_{\pi(g)}
\left[\prod_{g \in G'} p_{\pi(g)}\right]^2
\right) + O_{d', |J|, |K|, M_0}\left(\frac{1}{F(M_0)}\right).$$
Inserting this and \eqref{pag} back into \eqref{post-cauchy}, we can estimate \eqref{post-cauchy} by
$$ O_{M_d, d', |J|, |K|}\left( F(M_{d'})^{-1/2} \prod_{g \in G_{\subsetneq g_0}} p_{\pi(g)} \prod_{g \in G'} p_{\pi(g)} \right)
+ O_{d', |J|, |K|, M_0}\left(\frac{1}{F(M_0)}\right) .$$
Re-inserting those elements $g$ of $G$ for which $|g| = d'$ using \eqref{pe-big}, we can estimate this by
$$ O_{M_d, d', |J|, |K|}( F(M_{d'})^{-1/4} \prod_{g \in G} p_{\pi(g)})
+ O_{d', |J|, |K|, M_0}\left(\frac{1}{F(M_0)}\right)$$
(for instance).  By choosing $F$ sufficiently rapid depending on $d'$, $|J|$, $|K|$, we can write this as
$$ o_{M_d \to \infty; d', |J|, |K|}(\prod_{g \in G} p_{\pi(g)}) + O_{d', |J|, |K|, M_0}\left(\frac{1}{F(M_0)}\right).$$
Combining this with the bounds \eqref{contrib-1}, \eqref{fm0} we obtain \eqref{vk-count}, which closes the induction.
This completes the proof of Lemma \ref{gencount}, and hence Lemma \ref{count-lemma}.
\end{proof}


\begin{thebibliography}{10}

\bibitem{AS} M. Ajtai, E. Szemer\'edi, \emph{Sets of lattice points that form no squares}, Studia Scientarium Mathematicarum Hungarica. \textbf{9} (1974), 9--11.

\bibitem{chung}
F. Chung, \emph{Regularity lemmas for hypergraphs and quasi-randomness}, Random Struct. Alg. \textbf{2} (1991), 241--252.

\bibitem{efr}
P. Erd\"os, P. Frankl, V. R\"odl, \emph{The asymptotic number of graphs not containing a fixed subgraph and a problem for hypergraphs having no exponent}, Graphs Combin. \textbf{2} (1986), no. 2, 113--121.

\bibitem{frankl}
P. Frankl, V. R\"odl, \emph{The uniformity lemma for hypergraphs}, Graphs Combinat. \textbf{8}(4) (1992), 309--312. 

\bibitem{frankl02}
P. Frankl, V. R\"odl, \emph{Extremal problems on set systems}, Random Struct. Algorithms \textbf{20} (2002), no. 2, 131-164. 

\bibitem{furst}
H. Furstenberg, \emph{Ergodic behavior of diagonal measures and a theorem of Szemer\'edi on arithmetic progressions}, J. Analyse Math. \textbf{31} (1977), 204--256.

\bibitem{fk}
H. Furstenberg, Y. Katznelson, \emph{An ergodic Szemer\'edi theorem for commuting transformations}. J. Analyse Math. \textbf{34} (1978), 275--291.

\bibitem{fk2}
H. Furstenberg, Y. Katznelson, \emph{An ergodic Szemer\'edi theorem for IP-systems and combinatorial theory}, J. Analyse Math. \textbf{45} (1985), 117--168.

\bibitem{gowers-sz}
T. Gowers, \emph{Lower bounds of tower type for Szemer\'edi's uniformity lemma}, Geom. Func. Anal. \textbf{7} (1997), 322--337.

\bibitem{gowers}
T. Gowers, \emph{A new proof of Szemeredi's theorem}, GAFA \textbf{11} (2001), 465--588.

\bibitem{gowers-hyper}
T. Gowers, \emph{Hypergraph regularity and the multidimensional Szemer\'edi theorem}, preprint.

\bibitem{gt-primes}
B. Green, T. Tao, \emph{The primes contain arbitrarily long arithmetic progressions}, preprint.

\bibitem{krs-hyper}
Y. Kohayakawa, V. R\"odl, J. Skokan,
\emph{Hypergraphs, quasi-randomness, and conditions for regularity},
J. Combin. Theory Ser. A \textbf{97} (2002), no. 2, 307--352.

\bibitem{komlos}
J. Koml\'os, M. Simonovits, 
\emph{Szemer\'edi's regularity lemma and its applications in graph theory},
Combinatorics, Paul Erd\"os is eighty, Vol. 2 (Keszthely, 1993), 295--352, 
Bolyai Soc. Math. Stud., 2, J\'anos Bolyai Math. Soc., Budapest, 1996.

\bibitem{nrs}
B. Nagle, V. R\"odl, M. Schacht, \emph{The counting lemma for regular $k$-uniform hypergraphs}, to appear, Random Structures and Algorithms.

\bibitem{rodl-icm} 
V. R\"odl, \emph{Some developments in Ramsey theory}, Proceedings of the International Congress of Mathematicians, Vol. I, II (Kyoto, 1990), 1455--1466, Math. Soc. Japan, Tokyo, 1991.

\bibitem{rs}
V. R\"odl, M. Schacht, \emph{Regular partitions of hypergraphs}, preprint.

\bibitem{rstt}
V. R\"odl, M. Schacht, E. Tengan, N. Tokushige, \emph{Density theorems and extremal hypergraph problems}, preprint.

\bibitem{rodl}
V. R\"odl, J. Skokan, \emph{Regularity lemma for $k$-uniform hypergraphs}, to appear, Random Structures and Algorithms.

\bibitem{rodl2}
V. R\"odl, J. Skokan, \emph{Applications of the regularity lemma for uniform hypergraphs}, preprint.

\bibitem{roth}
K.F. Roth, \emph{On certain sets of integers}, J. London Math. Soc. \textbf{28} (1953), 245-252.

\bibitem{rsz}
I. Ruzsa, E. Szemer\'edi, \emph{Triple systems with no six points carrying three triangles}, Colloq. Math. Soc. J. Bolyai \textbf{18} (1978), 939--945.

\bibitem{shkredov}
I.D. Shkredov, \emph{On a problem of Gowers}, preprint.

\bibitem{soly-roth}
J. Solymosi, \emph{Note on a generalization of Roth's theorem}, Discrete and computational geometry, 825--827, Algorithms Combin. \textbf{25}, Springer Verlag, 2003.

\bibitem{soly-2}
J. Solymosi, \emph{A note on a question of Erdos and Graham}, Combinatorics, Probability and Computing \textbf{13} (2004), 263--267.

\bibitem{szemeredi}
E. Szemer\'edi, \emph{On sets of integers containing no $k$ elements in arithmetic progression},
Acta Arith. \textbf{27} (1975), 299--345.

\bibitem{szemeredi-reg}
E. Szemer\'edi, \emph{Regular partitions of graphs}, in ``Problem\'es Combinatoires et Th\'eorie des Graphes, Proc. Colloque Inter. CNRS,'' (Bermond, Fournier, Las Vergnas, Sotteau, eds.), CNRS Paris, 1978, 399--401.

\bibitem{tao:regularity}
T. Tao, \emph{Szemer\'edi's regularity lemma revisited}, preprint.

\bibitem{tao-multiprime}
T. Tao, \emph{The Gaussian primes contain arbitrarily shaped constellations}, preprint.

\bibitem{vdw}
B.L. Van der Waerden, \emph{Beweis einer Baudetschen Vermutung}, Nieuw. Arch. Wisk. \textbf{15} (1927), 212--216.

\bibitem{var}
P. Varnavides, \emph{On certain sets of positive density}, J. London Math. Soc. \textbf{39} (1959), 358--360.

\bibitem{van}
V. Vu, \emph{On a question of Gowers}, Ann. of Combinatorics, \textbf{6} (2002), 229--233.

\end{thebibliography}
\end{document}